\documentclass[a4paper]{article}

\usepackage{textcomp}

\usepackage{amsmath, amssymb, amsfonts}
\usepackage[english]{babel}
\def\<{\langle}\def\>{\rangle}
\def\R
{\mathbb{R}}\def\E
{\mathbb{E}}

\def\d{\frac{1}{2}\, }
\def \tr{\mathrm{\, trace}\, }
\title{Gaussian approximation of  Gaussian scale mixtures} 
\author{G\'erard Letac\thanks{Universit\'e Paul Sabatier, Toulouse} , H\'el\`ene Massam  \thanks{York University, Toronto}   }
   \date{}

\begin{document}

\maketitle
\begin{abstract}

For a given positive random variable $V>0$ and a given $Z\sim N(0,1)$ independent of $V$, we compute the scalar $t_0$ such that the distance   in the $L^2(\R)$ sense  between $Z V^{1/2}$ and $Z\sqrt{t_0}$ is minimal. We also consider the same problem in several dimensions when $V$ is a random positive definite matrix.

\noindent \textsc{Keywords:} Normal approximation,  Gaussian scale mixture, Plancherel theorem. 
\newline \textsc{ AMS Classification MSC2010:} Primary  62H17, Secondary 62H10.

\end{abstract}
\vspace{5mm}

{\it \noindent This work is dedicated to the memory of Frantisek Mat\'u$\check{s}$.}

\section{Introduction} Let $Z\sim N(0,I_n)$ be a standard Gaussian random variable in $\R^n$. Consider an independent random  positive definite matrix $V$ of order $n$ with distribution $\mu$. We call the distribution of $ V^{1/2}Z$ a Gaussian scale mixture, where $V^{1/2}$ is the unique positive definite matrix such that  $(V^{1/2})^2=V.$ Denote by $f$ the  density  of $ V^{1/2}Z$  in $\R^n$. In many practical circumstances, $\mu$ is not very well known, and $f$ is complicated. On the other hand, for $n=1$,  and
\begin{equation}\label{GENERIC}f(x)=\int_{0}^{\infty}e^{-\frac{x^2}{2v}}\frac{\mu(dv)}{\sqrt {2\pi v}}\end{equation}
we note that, as the logarithm of a Laplace transform,  $\log f(\sqrt x)$ is convex and thus the histogram of the symmetric density \eqref{GENERIC} looks
 like that  of a normal distribution. The central aim of the present paper is to say something  of the best normal approximation $N(0,t_0)$ of $f$ in the sense of $L^2(\R^n).$

In Section 2, we  recall some known facts and examples about the pair $(f,\mu)$ when $n=1.$ In Section 3,  our main result, for $n=1,$ is Theorem 3.1 in which  we show the existence of $t_0$,  its uniqueness and the fact that $t_0<\E(V)$. This theorem  also gives the equation, see \eqref{SOLVE},  that has to be solved to obtain  $t_0$ when $\mu$ is known. In Section 4 we consider the case $n\geq 2$ and investigate  the fact that several distributions of the random positive definite matrix $V$ can give the same Gaussian mixture $V^{1/2}Z.$  In Section 5, we  consider the problem of the Gaussian approximation of a Gaussian mixture in the more difficult case $n\geq 2.$  In that case, $t_0$ is a positive definite matrix, and in Theorem 5.2, we show the existence of $t_0$. Proposition 5.3 considers the particular case where $V$ is concentrated on the multiples of $I_n. $
A basic tool  we use in this paper is the Plancherel identity.

\section{The unidimensional case: a review}

A probability density $f$ on $\R$ is called a discrete Gaussian scale mixture if there  exist numbers $0<v_1<\cdots<v_n$ and $p_1,\ldots,p_n>0$ such that $p_1+\cdots+p_n=1$ and
$$f(x)=\sum_{i=1}^np_i\frac{1}{\sqrt{2\pi v_i}}e^{-\frac{x^2}{2v_i}}.$$
It easy to see that if 
$V\sim \sum_{i=1}^np_i\delta_{v_i}$ is independent of $Z\sim N(0,1)$ then the density of $Z V^{1/2}$ is $f.$ A way to see this is to observe that for all $s\in \R$ we have
$$\int_{-\infty}^{\infty}e^{sx}f(x)dx= \sum_{i=1}^np_ie^{\frac{s^2}{2}v_i}=\E(\E(e^{sZ V^{1/2}}|V))=\E(e^{sZ V^{1/2}}).$$ More generally, we will say that  the density $f$ is a \textit{Gaussian  scale  mixture} if there exists a probability distribution $\mu(dv)$ on $(0,\infty)$ such that \eqref{GENERIC} holds. 
As in the finite mixture case, if $V\sim \mu$ is independent of $Z\sim N(0,1)$ the density of $Z V^{1/2}$ is $f.$ To see this denote
\begin{equation}\label{LAPLACE}L_V(u)=\int_0^{\infty}e^{-uv}\mu(dv).\end{equation} Then
\begin{equation}\label{UNICITY}\int_{-\infty}^{\infty}e^{sx}f(x)dx=L_V(-s^2/2)=\E(e^{sZ V^{1/2}}).\end{equation}
For instance if $a>0$ and if \begin{equation}\label{FIRSTLAPLACE}f(x)=\frac{a}{2}e^{-a|x|}\end{equation} is the double  exponential density, then 
for $|s|<a$ we have $$\int_{-\infty}^{\infty}e^{sx}f(x)dx=\frac{a^2}{a^2-s^2}=L_V(-s^2/2)$$ where 
$$L_V(u)=\frac{a^2}{a^2+2u}=\frac{a^2}{2}\int_{0}^{\infty}e^{-vu-\frac{a^2}{2}v}dv.$$ This means that the mixing measure $\mu(dv)$ is an exponential distribution with mean $2/a^2. $

There are other  examples of pairs $(f,\mu)\sim (Z V^{1/2},V)$  in the literature. For instance, Palmer, Kreutz-Delgado and Makeig (2011) offer an interesting list of univariate mixing measures, containing also  some examples with $n>1.$  Another such list can be found in Gneiting (1997). Note that if $f$ is known then the distribution of $\log Z^2+\log V$ is known and finding the distribution $\mu$ or the distribution of $\log V$ is a problem of deconvolution. If its solution exists, it is unique, as shown for instance by \eqref{UNICITY}.

An example  of such a deconvolution is given by  West (1987), who extends \eqref{FIRSTLAPLACE} to $f(x)=Ce^{-a|x|^{2\alpha}}$ where  $0<\alpha<1$  as follows: he  recalls that for $A>0$ and $0<\alpha<1$, see Feller 1966, p. 424, there exists a probability density $g$, called a positive stable law, such that, for $\theta>0,$
\begin{equation}
\label{g}
\int_0^{\infty}e^{-t\theta}g(t)dt=e^{-A\theta^{\alpha}}.
\end{equation} 
If  in the equality above  we make the change of variable $t\to v=1/t$, let $\theta=x^2/2$ and define
 $\mu(dv)=C\sqrt{2\pi}g(1/v)v^{-3/2}dv$, where $C$ is such that $\mu(dv)$ is a probability,
we obtain
\begin{equation}
\label{int}
\int_{0}^{\infty}e^{-\d \frac{x^2}{v}}\frac{1}{\sqrt{2\pi v}}\mu(dv)=C e^{-2^{-\alpha}A|x|^{2\alpha}}.
\end{equation} 
Integrating both sides of \eqref{int} with respect to $x$ from $-\infty$ to $+\infty$, we obtain
$$C=\alpha\frac{A^{\frac{1}{2\alpha}}}{\sqrt{2}}\frac{1}{\Gamma(\frac{1}{2\alpha})}.$$

If $V\sim \mu$, its Laplace transform $L_V$ cannot be computed except for $\alpha=1/2$. For  $\alpha=1/2$ and $A$ arbitrary, one can verify that \eqref{g} is satisfied for 
	$$g(t)=\frac{A}{2\sqrt{\pi}}t^{-3/2}e^{-\frac{A^2}{4t}}.$$
	Then 
	$$\mu(dv)=\frac{A^2}{4}e^{-\frac{A^2}{4}v}{\bf 1}_{(0, +\infty)}(v) dv,$$ that is the mixing distribution is  an exponential distribution again.

Another elegant example of deconvolution is given by Stefanski  (1990) and Monahan and Stefanski (1992)
with the logistic distribution 
\begin{equation}\label{LOGISTIC}
f(x)=\frac{e^{x}}{(1+e^x)^2}=\sum_{n=1}^{\infty}(-1)^{n+1}ne^{-n|x|}.
\end{equation}
Using the representation of \eqref{FIRSTLAPLACE} as an exponential mixture of scale Gaussians, i.e.
\begin{equation*}
\frac{a}{2}e^{-a|x|}=\int_0^{+\infty}\frac{e^{-\frac{x^2}{2v}}}{\sqrt{2\pi v}}\frac{a^2}{2}e^{-\frac{a^2v}{2}}dv
\end{equation*}
 and applying it to $a=n$ in \eqref{LOGISTIC} above, we obtain
 \begin{eqnarray}
 f(x)=\sum_{n=1}^{\infty}(-1)^{n+1}n^2\int_0^{+\infty}\frac{e^{-\frac{x^2}{2v}}}{\sqrt{2\pi v}}e^{-\frac{n^2v}{2}}dv
 \end{eqnarray}
and thus,  if $\mu$ exists here, it must be 
\begin{equation}
\label{KS}
\mu(dv)=\left(\sum_{n=1}^{\infty}(-1)^{n+1}n^2e^{-\frac{n^2}{2}v}\right){\bf 1}_{(0, +\infty)}(v)dv
\end{equation}
which indeed exists since this is the Kolmogorov (1933) distribution, also called Kolmogorov-Smirnov distribution.  A direct proof that \eqref{KS} defines a probability  on $(0, +\infty)$  relies on  the following Jacobi  formula (see Hardy and Wright, 1938):
\begin{equation}
\label{ELLIPT}
\prod_{n=1}^{\infty}(1-q^{2n-1})^2(1-q^{2n})=\sum_{n=-\infty}^{\infty }(-1)^nq^{n^2}.
\end{equation} 
Taking $q=e^{-x/2}$, \eqref{ELLIPT} yields
\begin{equation*}
\prod_{n=1}^{\infty}(1-e^{-(2n-1)x/2})^2(1-e^{-nx})=\sum_{n=-\infty}^{\infty }(-1)^ne^{-n^2x/2}:=F(x).
\end{equation*} 
We observe that $F(0)=0$,  $F(+\infty)=1$, and $F$ is increasing as the  product of  increasing positive factors. Moreover,
$$F'(x)=-\frac{1}{2}\sum_{n-\infty}^{+\infty}(-1)^nn^2e^{-n^2x?2}=\sum_{n=1}^{+\infty}(-1)^{n+1}n^2e^{-n^2x/2}$$
is the density of \eqref{KS}.

\section {The normal approximation to the Gaussian scale mixture} The mixture $f$ as defined in \eqref{GENERIC} keeps some characteristics of the normal distribution: It is a symmetric density, $f(x)=e^{-\kappa(\frac{x^2}{2})}$ where $u\mapsto \kappa(u)$ is convex since 
$$e^{-\kappa(u)}=\int_0^{\infty}e^{-u/v}\frac{\mu(dv)}{\sqrt{2\pi v}}=\int_0^{\infty}e^{-uw}\nu(dw)$$ is the Laplace transform of the positive measure $\nu(dw)$ defined as the image of $\frac{\mu(dv)}{\sqrt{2\pi v}}$ by the map $u\mapsto w=1/v.$ 

As said in the introduction, in some practical applications, the distribution of $V$ is not very well known, and it is interesting to replace $f$ by the density of an ordinary normal distribution $N(0, t_0).$
The $L^2(\R)$ distance is well adapted to this problem. See Letac, Massam and Mohammadi (2018) for an example of the utilisation of this idea. We are going to prove the following result. 

\vspace{4mm}\noindent\textbf{Theorem 3.1.} If $f$ is defined by \eqref{GENERIC}, then 
\begin{enumerate}\item $f\in L^2(\R)$ if and only if $$\E\left(\frac{1}{\sqrt{V+V_1}}\right)<\infty$$ when $V$ and $V_1$ are independent with the same distribution $\mu.$
\item If $f\in L^2(\R)$, there exists a unique $t_0=t_0(\mu)>0$ which minimizes 
$$t\mapsto I_V(t)=\int_{-\infty}^{\infty}\left[f(x)-\frac{1}{\sqrt{2\pi t}}e^{-\frac{x^2}{2t}}\right]^2dx.$$ 
\item The scalar $y_0=1/t_0$ the unique positive solution  of the equation 
\begin{equation}\label{SOLVE}\int_0^{\infty}\frac{\mu(dv)}{(1+vy)^{3/2}}=\frac{1}{2^{3/2}}.
\end{equation} 
In particular, if $\mu_{\lambda} $ is the distribution of $\lambda V$, then       $t'=t_0(\mu_{\lambda})=\lambda t_0(\mu).$

\item The value of $I_V(t_0)$ is $$I_V(t_0)=\sqrt{\frac{2}{\pi}}\left(\E\left(\frac{1}{\sqrt{V+V_1}}\right)-2\E\left(\frac{1}{\sqrt{V+t_0}}\right)+\frac{1}{\sqrt{2t_0}}\right)$$and
\begin{equation}\label{LAMBDA} I_{\lambda V}(t')=\frac{1}{\sqrt{\lambda}}I_V(t_0).\end{equation}

\item Finally
$t_0\leq \E(V).$\end{enumerate}

\vspace{4mm}\noindent\textbf{Proof.}  Recall that if $g\in L^2(\R)\cap L^1(\R)$ and if $\hat{g}(s)=\int_{-\infty}^{\infty}e^{isx}g(x)dx$, then  Plancherel theorem says that
\begin{equation}\label{PLANCHEREL}\frac{1}{2\pi}\int_{-\infty}^{\infty}|\hat{g}(s)|^2ds=\int_{-\infty}^{\infty}|g(x)|^2dx.\end{equation}
 Furthermore if $g\in L^1(\R),$ then $g\in L^2(\R)$ if and only if $\hat{g}\in L^2(\R).$

Let us apply \eqref{PLANCHEREL} first to $g=f.$ From \eqref{GENERIC} and \eqref{UNICITY}, we have $\hat{f}(s)=L_V(s^2/2)$. Then

\begin{eqnarray*}
\int_{-\infty}^{\infty}\hat{f}^2(s)ds&=&\int_{-\infty}^{\infty}L_V^2(s^2/2)ds=\sqrt{2}\int_{0}^{\infty}L(u)^2\frac{du}{\sqrt{u}}\\&=&\sqrt{2}\int_{0}^{\infty}\E(e^{-u(V+V_1)})\frac{du}{\sqrt{u}}=\sqrt{2\pi}\, \E(\frac{1}{\sqrt{V+V_1}})\end{eqnarray*}
where the last equality is obtained by recalling that $\int_0^{+\infty} e^{-uv}\frac{dv}{\sqrt{v}}=\frac{\sqrt{\pi}}{\sqrt{u}}.$
Thus  statement 1.  of the theorem is proved.

To prove 2., 3.  and 4.,
we  apply \eqref{PLANCHEREL} to $g(x)=f(x)-\frac{1}{\sqrt{2\pi t}}e^{-\frac{x^2}{2t}}$ for which $\hat{g}(s)=L(s^2/2)-e^{-ts^2/2}.$ As a consequence
$$I_V(t)=\frac{1}{2\pi}\int_{-\infty}^{\infty}\left[L_V(s^2/2)-e^{-ts^2/2}\right]^2ds=\frac{1}{\pi}\int_{0}^{\infty}\left[L_V(u)-e^{-tu}\right]^2\frac{du}{\sqrt{2u}}$$
and
\begin{equation}
\label{i'vt}
I_V'(t)=\frac{\sqrt{2}}{\pi}\int_{0}^{\infty}\left[L_V(u)-e^{-tu}\right]e^{-tu}\sqrt{u}du. 
\end{equation} 
Since $\int_{0}^{\infty}e^{-2tu}\sqrt{u}du=\frac{\Gamma(3/2)}{(2t)^{3/2}}$ and since 
$$\int_{0}^{\infty}L_V(u)e^{-tu}\sqrt{u}\;du=\int_{0}^{\infty}\int_{0}^{\infty}
e^{-u(v+t)}\sqrt{u}\;du\, \mu(dv)=\Gamma(3/2)\int_{0}^{\infty}\frac{\mu(dv)}{(t+v)^{3/2}},$$ then $I_V'(t)=0$ if and only if 
$\int_{0}^{\infty}\frac{\mu(dv)}{(t+v)^{3/2}}=\frac{1}{(2t)^{3/2}}.$
We can rewrite this equation in $t$ as $F(1/t)=1/2^{3/2}$ where 
$F(y)=\int_0^{\infty}\frac{\mu(dv)}{(1+vy)^{3/2}}.$
Thus \eqref{i'vt} can be rewritten
\begin{equation}
\label{re:iv't}
I_V'(t)=\frac{\sqrt{2}}{\pi}\frac{\Gamma(3/2)}{t^{3/2}}\left[F\left(\frac{1}{t}\right)-\frac{1}{2^{3/2}}\right].
\end{equation}
Since $0<1/2^{3/2}<1,$  $F(0)=1,$ $ \lim_{y\to \infty}F(y)=0$ and 
$$F'(y)=-\frac{3}{2}\int_0^{\infty}\frac{v\mu(dv)}{(1+vy)^{5/2}}<0,$$
it follows that $I_V'$ has only one zero $t_0$ on $(0,\infty)$ and from \eqref{re:iv't}, it is easy to see from the sign of $I_V'$ that $I_V$ reaches its minimum at $t_0.$

 To show 5., we will apply Jensen inequality $f(\E(X))\leq \E(f(X))$  to the convex function $f(x)=x^{-3/2}$ and to  the random variable $X=1+y_0V$. From
$$\frac{1}{(1+y_0\E(V))^{3/2}}\leq \E\left(\frac{1}{(1+y_0V)^{3/2}}\right)=\frac{1}
{2^{3/2}}$$ it follows that  $2\leq 1+y_0\E(V)$ and $t_0=1/y_0<\E(V).$  \hfill $\square$

\vspace{4mm}\noindent\textbf{Example 1.} Suppose that 
$\Pr(V=1)=\Pr(V=2)=1/2.$ Let us compute $t_0$ and $I(t_0).$ With the help of Mathematica, we see that
the solution of $$\frac{1}{2(1+t)^{3/2}}+\frac{1}{2(2+t)^{3/2}}=\frac{1}{(2t)^{3/2}}$$ is  $t_0=1.39277$. Finally 
$$I_V(t_0)=\sqrt{\frac{2}{\pi}}\left(\frac{1}{4\sqrt{2}}+\frac{1}{2\sqrt{3}}+\frac{1}{8}-\frac{1}{\sqrt{1+t_0}}-\frac{1}{\sqrt{2+t_0}}+\frac{1}{\sqrt{2t_0}}\right)=0.00019,$$
which is very small. 

\vspace{4mm}\noindent\textbf{Example 2.} Suppose that $V$ is uniform on $(0,1)$ Then $$t_0=0.36678, \ I_V(t_0)=0.0182.$$ If $V$ is uniform on $[0,a]$, then from Part 4 of  Proposition  3.1, we have  $t_0=a\times 0.36678.$

\vspace{4mm}\noindent\textbf{Example 3.} If $V$ follows the standard exponential distribution with density $f(v)=e^{-v}{\bf 1}_{(0,+\infty)}(v)$, then $$t_0=0.524, \ I_V(t_0)=0.0207.$$

\section{Scale mixtures in the Euclidean case and non  identifiability} 
Denote by $\mathcal{S}$ the linear space of symmetric real matrices of dimension $n$ equipped with the scalar product $\<s,s_1\>=\tr (s s_1)$ 
 and by $\mathcal{P}$ the convex cone of real positive definite matrices of order $n$. Thus the norm of $s$ is $||s||=\sqrt{\tr s^2}.$ We denote by $dv$ the Lebesgue measure on  $\mathcal{S}$ associated to its Euclidean structure, namely such that the mass of a unit cube is one.

We use the symbol $a^*$ for the transposed  matrix of any matrix $a.$ As said before, if $v\in \mathcal{P}$ we denote by $v^{1/2}$ the unique element of $\mathcal{P}$ whose square is $v.$ 
 It is sometimes considered that any non singular matrix $a$ such that $v=aa^*$ should be called a generalized square root of $v.$
The Cholesky decomposition $v=tt^*$ of $v$ into a product of a upper triangular matrix $t$ with positive coefficients on the diagonal with its transposed matrix $t^*$ offers an example of such a generalized square root. It can be remarked that in practice  the calculation of $t$ is easier than the calculation of $v^{1/2}$. We denote by $\mathbb{O}(n)$ the orthogonal group of $n\times n$ matrices $u$ such that $u^*u=I_n.$

In this section we  define the scale mixtures  of the standard normal  distribution in $\R^n$ and we observe the phenomena of non identifiability: that is,  different  distributions of $V$ can give the same mixture.

 \subsection{Scale mixtures of the normal distribution in $\R^n$.}

 A scaled Gaussian mixture $f$ on $\R^n$  is the density of  a random variable $X$ on $\R^n$ of the form $X=V^{1/2}Z$ where $V\sim \mu$ is a random matrix in $\mathcal{P}$ independent of the standard random Gaussian variable $Z\sim N(0,I_n).$  
We give some properties of such a mixture in the following proposition.

\vspace{4mm}\noindent\textbf{Proposition 4.1.} Let $A$ be a random nonsingular square matrix of order $n$, independent of $Z\in \R^n\setminus \{0\}$   and such that $uZ\sim Z$ for all $u\in \mathbb{O}(n).$ Let $V=AA^*$. Then the following holds.

\begin{enumerate}

\item $AZ\sim V^{1/2}Z$, that is,  if we replace $V^{1/2}$ by any generalized square root $A$ of $V$, the distribution of $AZ$ remains the same.

\item If $AZ\sim Z$ then $ \Pr(V=I_n)=1.$  In other terms,  $AZ\sim Z$ if and only if  $ \Pr(AA^*=I_n)=1$, i.e
$A\in \mathbb{O}(n)$ almost surely. 

\end{enumerate}

\vspace{4mm}\noindent\textbf{Proof.} To prove 1.,  observe that $U=V^{-1/2}A$ is in the orthogonal group $\mathbb{O}(n)$. Let $\mu(dv)K(v,du)\nu(dz)$ denote the joint distribution of $(V,U,Z)$.

Then if $h$ is a bounded function on $\R^n$, 
\begin{eqnarray}
\E(h(AZ))&=&\E(h(V^{1/2}UZ))\nonumber\\&=&\int_{\mathcal P} \mu(dv)\int_{\mathbb{O}(n)} K(v,du)\int_{\R^n}h(v^{1/2}uz)\nu(dz)\nonumber\\
&=&\int_{\mathcal P} \mu(dv)\int_{\mathbb{O}(n)} K(v,du)\int_{\R^n}h(v^{1/2}z_1)\nu(dz_1)\label{z1}\\
&=&\int_{\mathcal P} \mu(dv)\int_{\R^n}h(v^{1/2}z_1)\nu(dz_1)=\E(h(V^{1/2}Z)),\label{u}
\end{eqnarray}
where in \eqref{z1}, $z_1=uz$, and  \eqref{u} follows from  $\int_{\mathbb{O}(n)} K(v,du)=1.$

To prove 2., consider also $\varphi(s)=\mathbb\E(e^{i\<s,Z\>}).$ Since  $uZ\sim Z$ for all $u\in \mathbb{O}(n)$ there exists a real function $g$ defined on $[0,\infty)$ such that 
$\varphi(s)=g(\|s\|^2).$ Since $Z\sim AZ$ we can write 
\begin{equation}\label{GV}
g(\|s\|^2)= \E(g(s^*Vs))\;.
\end{equation}

Next, let us show that if $R\geq 0$ is independent of 
$Z=(Z_1,\ldots,Z_n)$ and if $Z_1R\sim Z_1$ then $\Pr(R=1)=1$. Indeed, for $t\geq 0$  we have that $\E(|Z_1|^{it})=\E(|Z_1|^{it})\E(R^{it})$. Since there exists $0<t_0\leq \infty$ such that $\E(|Z_1|^{it})\neq 0$ for $0\leq t<t_0$, it holds that 
$\E(R^{it})=1$ for $0\leq t<t_0$. This implies  that $\Pr(R>0)=1$ and $ 0=1-\Re(\E(R^{it}))=\E(1-\cos (t\log R))$ or $\Pr(t\log R\in 2\pi \mathbb{Z})=1$ for $0\leq t<t_0.$ We deduce easily  that $\Pr(R=1)=1.$

Now denote $V=(V_{ij})_{1\leq i,j\leq n}$ and apply the above observation to $R=\sqrt{V_{11}}$ by taking $s=(t,0,\ldots,0)$ in \eqref{GV}. We obtain

$$\E(e^{it Z_1})=\varphi((t,0,\ldots,0))=g(t^2)=\E(g(t^2 V_{11}))=\E(e^{it \sqrt{V_{11}}Z_1})$$ which implies $Z_1\sim {V_{11}}Z_1$ and $\Pr(V_{11}=1)=1.$
Similarly $\Pr(V_{ii}=1)=1$ for all $i=2,\cdots,n.$

Finally, we consider $R=\sqrt{1+V_{12}}$ and we take   $s=(t/\sqrt{2},t/\sqrt{2},\ldots,0)$ in \eqref{GV}. Using the fact that $(Z_1+Z_2)/\sqrt{2}\sim Z_1$ we write 
\begin{eqnarray*}\E(e^{it Z_1})&=&\E(e^{it (Z_1+Z_2)/\sqrt{2}})=\varphi((t/\sqrt{2},t/\sqrt{2},\ldots,0))\\&=&\E(g(  \d t^2 (V_{11}+V_{22}+2V_{12}))=\E(g(   t^2 (1+V_{12}))\\
&=&\E(e^{itZ_1\sqrt{1+V_{12}}})\end{eqnarray*} and we get $\Pr(V_{12}=0)=1.$ Similarly $\Pr(V_{ij}=0)
=1$ for $i\neq j$ and finally $\Pr(V=I_n)=1$ as desired. \hfill $\square$

\subsection{Nonidentifiability}

In Example 4 below, we show that for $ n\geq 2,$  the measure $\mu$ which generates a given $f$  as in \eqref{GENERIC} may not be unique. Theorem 4.2 gives a more general result. We denote  by $\omega$ the uniform probability, or Haar probability,  on  $\mathbb{O}(n)$ and by $\mathcal{D}$ the set of diagonal matrices 
$b=\mathrm{diag}(b_1,\ldots,b_n)$ such that $0<b_1\leq b_2\leq \ldots\leq b_n.$  
It is a well known fact  that  if $V=U^*BU$ with $U\in  \mathbb{O}(n)$ and $B\in \mathcal{D}$ then  $u^*Vu\sim V$ for all $u\in \mathbb{O}(n)$  if and only if $U\sim \omega$ and $B$ are independent (in this case, the distribution of $V$ is determined by the distribution  of its set of eigenvalues determined by $B$). While the 'if' part is clear, a short proof of the 'only if ' part is as follows: consider $\alpha(db)K(b,du)\sim(B,U)$
and $\mu\sim V.$ For any $h$ bounded continuous on $\mathcal{P}$ and any $u_0\in \mathbb{O}(n)$ we write\begin{eqnarray*}\int_{\mathcal{P}}h(v)\mu(dv)&=&\int_{\mathcal{P}}h(u^*_0vu_0)\mu(dv)\\&=&\int_{\mathcal{D}}\left(\int_{\mathbb{O}(n)}h(u_0^*u^*buu_0)K(b,du)\right)\alpha(db)\\&=&\int_{\mathcal{D}}\left(\int_{\mathbb{O}(n)}h(u^*bu)K(b,d(uu_0^*)\right)\alpha(db)\end{eqnarray*}
This shows that, $\alpha$ almost surely, the probability $ K(b,du)$ on $\mathbb{O}(n)$ is invariant by $u\mapsto 
uu_0^*$ for all $u_0\in \mathbb{O}(n)$  and is equal to $\omega$ by uniqueness of the Haar probability on $ \mathbb{O}(n).$

Finally, for $a_1,\ldots a_n>0$ given, we recall the definition of the Dirichlet distribution $D(a_1,\ldots,a_n))$ of the variable $(X_1,\ldots,X_n)$ on the simplex $$T_n=\{(x_1,\ldots,x_n)\in (0,\infty)^n\ ;\ x_1+\cdots+x_n=1\}:$$  the density of $(X_2,\ldots,X_n)$ is proportional to $$(1-(x_2+\cdots+x_n)^{a_1-1}x_2^{a_2-1}\ldots x_n^{a_n-1}.$$

\vspace{4mm}\noindent\textbf{Theorem 4.2.} Suppose that a probability $\mu(dv)$ on $\mathcal{P}$ is invariant by the transformations  $v\mapsto uvu^*$ for any $u\in \mathbb{O}(n).$  Then we have the following.

\begin{enumerate}

\item  Let $V\sim \mu$. Then there exists  a unique probability $\nu_{\mu}(d\lambda)$ on $(0,\infty)$ such that if $\Lambda\sim \nu_{\mu}$  and if $V$ and $\Lambda$ are independent of $Z\sim N(0,I_n)$, then $$V^{1/2}Z=\Lambda^{1/2}Z.$$
\item In the special case where  $b=\mathrm{diag}(b_1.\ldots,b_n)\in\mathcal{D}$ is fixed let $\mu_b$ be the distribution  in $\mathcal{P}$ of $ U^*bU$ where $U\sim \omega $ . For   $(X_1,\ldots, X_n)\sim D(\d,\d,\ldots,\d)$,  denote by $\rho_b(d\lambda)$ the distribution of $b_1X_1+\cdots+b_nX_n$. Then \begin{equation}\label{SUPER}\rho_b=\nu_{\mu_b}.\end{equation}
\item If $\alpha(db)$ is a probability on $\mathcal{D}$, denote by $\mu$ the distribution of $V=U^*BU$ where $B\sim \alpha$ and $U\sim \omega$ are independent.  Then 
\begin{equation}\label{REP}\nu_{\mu}(d\lambda)=\int_{\mathcal{D}}\alpha(db)\rho_b(d\lambda).\end{equation}

\end{enumerate}

\vspace{4mm}\noindent\textbf{Proof.}  
We begin with a remark. 
Consider  the Fourier transform of $ V^{1/2}Z$
defined for $s\in \R^n$ by 
$\varphi(s)=\E(e^{is^* V^{1/2}Z})=\E(e^{-\d s^*Vs}).$
For $u\in \mathbb{O}(n)$ the fact that $u^*Vu\sim V$ implies that 
$\varphi(us)=\varphi(s).$ This implies in turn that $\varphi(s)$ is a function of $\|s\|$ only, or that there exists a function $L$ such that $\varphi(s)=L(\d\|s\|^2).$  Recall that we intent to show the existence of a positive random variable $\Lambda$ such that $L(\d\|s\|^2)=\E(e^{-\d\Lambda\|s\|^2})$ that is, that $L$ is a Laplace transform. Actually this point is not immediate, and we start the proof of the theorem by showing  \eqref{SUPER} first. 

Let $V=U^*bU$ with $U\sim \omega$ and consider the Fourier transform $\varphi(s)$ of $V^{1/2}Z$, namely
\begin{equation}\label{ORTHO}\varphi(s)=\E(e^{-\d (Us)^*bUs})=\E(e^{-\d(b_1(Us)_1^2+\cdots +b_n(Us)_n^2)})\end{equation} where $Us=((Us)_1,\ldots,(Us)_n).$ Now we observe that $(Us)/\|s\|$ is uniformly distributed on the unit sphere of $\R^n.$  If $Y=(Y_1,\ldots,Y_n)\sim N(0,I_n)$ then  $Y/\|Y\|$ is also uniformly distributed on the sphere and it is a classical fact that 
$$(X_1,\ldots,X_n)=\frac{(Y_1^2,\ldots, Y_n^2)}{Y_1^2+\cdots+Y_n^2}\sim D(\d,\ldots,\d)$$ Therefore 
$$\frac{1}{\|s\|^2}(Us)^*b(Us)\sim b_1X_1+\cdots+b_nX_n\sim \rho_b$$ and $\varphi(s)=\int_0^{\infty}e^{-\d \|s\|^2 \lambda}\rho_b(d\lambda)$, which is a reformulation of \eqref{SUPER}.  Note that in this particular case where $V=U^*bU$ then $L$ is the Laplace transform of $\rho_b.$

To prove 3.,   we simply condition by $B$ and use \eqref{SUPER} to obtain $$\varphi(s)=\E(e^{-\d (Us)^*B(Us)}) =\int_{\mathcal{D}}\left(\int_0^{\infty}e^{-\d \|s\|^2\lambda}\rho_b(d\lambda)\right)\alpha(db)$$ which proves \eqref{REP}. 

Recall that any random variable $V$ on $\mathcal{P}$ such that  $u^*Vu\sim V$ for all $u\in \mathbb{O}(n)$ has the above form $U^*BU$  where $B\sim \alpha(db)$ is random and independent of $U\sim \omega.$ This shows that 3. implies 1. 
\hfill $\square$
 
\vspace{4mm}\noindent\textbf{Corollary 4.3.} If $V\sim   uVu^*$ for any $u\in \mathbb{O}(n)$ and has distribution $\mu$ then the density $f$ of $V^{1/2}Z$ where $Z\sim N(0,I_n)$ is independent of $V$ has the form $f(x)=L_1(\|x\|^2/2).$  More specifically
\begin{equation}\label{REPDENSITY}f(x)=\int_0^{\infty}e^{-\frac{\|x\|^2}{2\lambda}}
\frac{\nu_{\mu}(d\lambda)}{\sqrt{2\pi \lambda}}.\end{equation}

\vspace{4mm}\noindent\textbf{Remarks.} 
\begin{enumerate}\item 
Note that in Corollary 4.3  the function $L_1$ is the Laplace transform of the image $m(dy)$ of the measure 
$\frac{\nu_{\mu}(d\lambda)}{\sqrt{2\pi \lambda}}$ by the map $\lambda\mapsto y= 1/2\lambda.$ Since  in general \eqref{REP}  is not easy to apply, 
this  offers, in some cases, a way to compute $\nu_{\mu}(d\lambda)$, when $f$  and $L_1$ are known, and when $m$ is obvious. Example 4 below will be obtained by this technique with $L_1(s)=(1+2s)^{-p}$ with $p>n/2.$

\item For $n\geq 3$ it is difficult to give the density of $\rho_b(d\lambda)$ explicitly. For $n=2$  it is the image of the beta distribution on $(0,1)$ with parameters (1/2,1/2) by the affinity $t\mapsto \lambda=(1-t)b_1+tb_2$ :
$$\rho_b(d\lambda)=\frac{1}{\pi\sqrt{(b_2-\lambda)(\lambda-b_1)}}1_{(b_1,b_2)}(\lambda)d\lambda.$$ For instance if $\alpha(db_1,db_2)=\alpha_1(db_1)K(b_1,db_2)$ is the joint distribution of $ B=\mathrm{diag}(B_1,B_2),$  formula \eqref{REP} implies 
$\nu_{\mu}(d\lambda)$ has density
$$\frac{1}{\pi}\int_0^{\lambda}\left(\int_{\lambda}^{\infty}\frac{K(b_1,db_2)}{\sqrt{b_2-\lambda}}\right)\frac{\alpha_1(db_1)}{\sqrt{\lambda-b_1}}.$$

 \item Another approach to  formula \eqref{SUPER} is possible using zonal polynomials.

Indeed for any symmetric matrices $a$ and $b$ of order $n$ we can write
$$\int_{\mathbb{O}(n)}e^{\tr u^*bua}\omega(du)=\sum_{\kappa}\frac{C_{\kappa}(a)C_{\kappa}(b)}{|\kappa|! C_{\kappa}(I_n)}.$$ Equality \eqref{ORTHO} suggests to apply this identity to the matrices $a=-ss^*/2$ and $b\in \mathcal{D}.$ Fortunately the zonal polynomials are simple when computed on $a,$ a matrix of rank one. More specifically $C_{\kappa}(a)=0$ except when $\kappa=(m,0,0,\ldots,0)$ where $m$ is a non negative integer. In this case,  by a reasoning similar  to that in the proof of \eqref{SUPER}, we have
$$ \frac{C_{\kappa}(a)}{|\kappa|! C_{\kappa}(I_n)}=\frac{(-1)^m}{2^mm!}\int_{\mathbb{O}(n)}
(us)_1^{2m}\omega(du)=\frac{(-1)^m\|s\|^{2m}}{2^mm!}\E(X_1^m)$$
where $X_1\sim \beta(\d,\d(n-1))$.
However, the computation of $$c_m(b_1,\ldots,b_n)=C_{(m,0,0,\ldots,0)}(\mathrm{diag}(b_1,\ldots,b_n))$$ is the real difficulty and using the Pochhammer symbol $(x)_n=\Gamma(n+x)/\Gamma(x),$  one can only write

$$\E(e^{-\d (Us)^*bUs})=\sum_{m=0}^{\infty}\frac{(-1)^m\|s\|^{2m}(1/2)_m}{2^m(n/2)_m m!}c_m(b_1,\ldots,b_n).$$

\item An interesting question is the following: suppose that more generally $V\sim \mu$ and $V_1\sim \mu_1$  in $\mathcal{P}$ are such that $V^{1/2}Z\sim V_1^{1/2}Z$ with $Z\sim N(0,I_n)$ independent of $V$ and $V_1.$ We do not assume here that $\mu$ and $\mu_1$ are invariant by $\mathbb{O}(n).$ Consider the Laplace transforms $L_{\mu}( a)=\int_{\mathcal{P}}e^{-\tr (av)}\mu(dv)$ and $L_{\mu_1}$ defined at least on the closed convex cone $\overline {\mathcal{P}}$ of the semi positive definite matrices of order $n$. Then $V^{1/2}Z\sim V_1^{1/2}Z$ implies that for any $s\in \R^n$ we have $$L_{\mu}(\d  ss^*)=L_{\mu_1}(\d  ss^*)$$ which means that $L_{\mu}$ and $L_{\mu_1}$ coincide on the matrices $a\in \overline {\mathcal{P}}$
 of rank one. As we have just seen in Theorem 4.2 it does not imply  $\mu=\mu_1.$
This raises the following problem: given $\mu$, describe the extreme points of  the convex set of probabilities $\mu_1$ such that $L_{\mu}$ and $L_{\mu_1}$ coincide on the matrices $a\in \overline {\mathcal{P}}$
 of rank one. 

\end{enumerate}

\subsection{An explicit example of non identifiability.} We will now give an example of two different measures $\mu_1$ and $\mu_2$ giving the same scale mixture of Gaussian variables. 

\vspace{4mm}\noindent\textbf{Example 4.} Let $p>n/2$  and consider the probability on $\R^n$ with density
\begin{equation}\label{ALEXA}
f(x)=\frac{C}{(1+||x||^2)^p},
\end{equation}
where $C$ will be computed below. Then consider  two probability measures $\mu_1$ and $\mu_2.$ The first is 
\begin{eqnarray}
\mu_1(dv)&=&\frac{(\det (v))^{-p+\frac{1}{2}-\frac{n+1}{2}}}{2^{n(p-\d) }\Gamma_{\mathcal{P}}(p-\frac{1}{2})}\exp\{ -\frac{1}{2}\tr(v^{-1}) \}{\bf 1}_{\mathcal{P}}dv, 
\end{eqnarray}
where $\Gamma_{\mathcal{P}}(t)=(2\pi)^{\d n(n-1)}\prod_{j=1}^d\Gamma(t-\frac{j-1}{2}).$ Therefore $V^{-1}$ follows a Wishart distribution with shape parameter $p-\d.$ 
The second is defined by  $\mu_2(dv)\sim {\Lambda}I_n$ where $\Lambda$ has density
$$\frac{\lambda^{-p+\frac{n}{2}-1}}{2^{p-\frac{n}{2}}\Gamma(p-\frac{n}{2})}e^{-\frac{1}{2\lambda}}{\bf 1}_{(0,+\infty)}(\lambda), $$i.e. $\Lambda^{-1}$ follows a Gamma distribution, with shape parameter  $p-\d n$.
For $x\in \R^n$   and  $i=1,2$, we now show that
\begin{equation}
\label{mixwish}
\int_{\mathcal{P}}\frac{e^{-\frac{x^*v^{-1}x}{2}}}{(2\pi)^{n/2}(\det v)^{1/2}}\mu_i(dv)=f(x)
\end{equation}where $f$ is defined by \eqref{ALEXA}.
For $i=1$, making the change of variable $y=v^{-1}$, the left-hand side of \eqref{mixwish} becomes
\begin{eqnarray}
&&\int_{\mathcal{P}}\frac{(\det y)^{1/2}e^{-\frac{x^*yx}{2}}}{(2\pi)^{n/2}} 
\frac{(\det (y))^{p-\frac{1}{2}-\frac{n+1}{2}}}{2^{n(p-\d) }\Gamma_{\mathcal{P}}(p-\frac{1}{2})}\exp\{ -\frac{1}{2}\tr y \}dy\nonumber\\
&=&\int_{\mathcal{P}}\frac{(\det (y))^{p-\frac{n+1}{2}}}{(2\pi)^{n/2}2^{n(p-\d) }\Gamma_{\mathcal{P}}(p-\frac{1}{2})}e^{-\frac{1}{2}\tr(y, I_n+xx^*)}dy\nonumber\\
&=&\frac{2^{np }\Gamma_{\mathcal{P}}(p)}{(2\pi)^{n/2}2^{n(p-\d) }\Gamma_{\mathcal{P}}(p-\frac{1}{2})}\det(I_n+xx^*)^{-p}\nonumber\\
&=&\frac{2^{np}\Gamma_{\mathcal{P}}(p)}{(2\pi)^{n/2}2^{n(p-\d) }\Gamma_{\mathcal{P}}(p-\frac{1}{2})}\frac{1}{(1+||x||^2)^p}=\frac{1}{(2\pi)^{n/2}}\frac{\Gamma(p)}{2^{-\frac{n}{2}}\Gamma(p-\frac{n}{2})}\frac{1}{(1+||x||^2)^p},\nonumber
\end{eqnarray}
yielding 
$C=\frac{1}{(2\pi)^{n/2}}\frac{\Gamma(p)}{2^{-\frac{n}{2}}\Gamma(p-\frac{n}{2})}.$
For $i=2$, making the change of variable $y=\frac{1}{\lambda}$, the left-hand side of \eqref{mixwish} becomes
 \begin{eqnarray}
&&\int_0^{+\infty}\frac{e^{-\frac{x^tx}{2\lambda}}}{(2\pi)^{n/2}\lambda^{\frac{n}{2}}} 
\frac{\lambda^{-p+\frac{n}{2}-1}}{2^{p-\frac{n}{2}}\Gamma(p-\frac{n}{2})}e^{-\frac{1}{2\lambda}}{\bf 1}_{(0,+\infty)}(\lambda)\nonumber\\
&=&\frac{1}{(2\pi)^{n/2}\Gamma(p-\frac{n}{2})}
\int_0^{+\infty}\frac{\lambda^{-p-1}}{2^{p-\frac{n}{2}}}e^{-\frac{1}{2\lambda}(1+||x||^2)}d\lambda\nonumber\\
&=&\frac{1}{(2\pi)^{n/2}\Gamma(p-\frac{n}{2})}\int_0^{+\infty}
\frac{y^{p-1}}{2^{p-\frac{n}{2}}}e^{-\frac{y}{2}(1+||x||^2)}dy\nonumber\\
&=&\frac{\Gamma(p)}{(2\pi)^{n/2}2^{-\frac{n}{2}}\Gamma(p-\frac{n}{2})}\frac{1}{(1+||x||^2)^p}
\end{eqnarray}
Therefore, with the notation of Theorem 4.2 we have proved that  if $\mu_2\sim \Lambda I_n$ then $\Lambda\sim \nu_{\mu_1}.$
\section{Existence of the best normal approximation in the Euclidean case}

In this section, we study  the conditions that the  distribution $\mu(dv)$ on $\mathcal{P}$ must satisfy to garantee that the density  $f$ of $V^{1/2}Z$ is  in  $L^2(\R^n)$ when $V\sim \mu$ and $Z\sim N(0,I_n)$ are independent. We also find a Gaussian law $N(0,t_0)$ on $\R^n$ which is the closest to $f$ in the $L^2(\R^n)$ sense. We consider also the particular case where $V^{1/2}Z=\Lambda^{1/2} Z$ where $\Lambda$ is a random scalar.

\subsection{Best approximation}

We first recall  two simple formulas. 

\vspace{4mm}\noindent\textbf{Lemma 5.1.} Let $A\in \mathcal{P}.$ Then
$$\int_{\R^n}e^{-\d s^*As}ds=\frac{(2\pi)^{n/2}}{\sqrt{\det A}},\ \ \ \int_{\R^n}e^{-\d s^*As}ss^*ds=\frac{(2\pi)^{n/2}}{\sqrt{\det A}}A^{-1}.$$ 

\vspace{4mm}\noindent\textbf{Proof.} Without loss of generality, we may assume that $A$ is diagonal, and the proof is obvious in this particular case.

We next state that there exists a matrix $v=t_0$ such that the $L^2$ distance between the multivariate Gaussian mixture $f(x)$ and the Gaussian distribution $N(0,t_0)$ is minimum.

\vspace{4mm}\noindent\textbf{Theorem 5.2.} Let $\mu(dv)$ be a probability distribution on the convex cone $\mathcal{P}$. Let
$f(x)$ denote the density of the random variable $X=V^{1/2}Z$ of $\R^n$ where $V\sim \mu$ is independent of $Z\sim N(0,I_n).$ Then \begin{enumerate}\item 
$f\in   L^2(\R^n)$ if and only if $\E\left(\frac {1}{\det\sqrt{V+V_1}}\right)<\infty$
where $V$ and $V_1$ are independent with the same distribution $\mu.$ 
\item For $f\in   L^2(\R^n),$ consider the function $I$ defined on $\mathcal{P}$ by
\begin{equation}
\label{ioft}
t\mapsto I(t)=\int_{\R^n}\left[f(x)-\frac{1}{\sqrt{(2\pi)^n\det t}}e^{-\d x^*t^{-1}x}\right]^2dx.
\end{equation}
Then $I$ reaches its minimum at some $t_0$, and this $t_0$ is  a solution in $\mathcal{P}$ of the following equation in $t\in \mathcal{P}:$
\begin{equation}\label{DIFFn}\int_{\mathcal{P}}\frac{(v+t)^{-1}}{\sqrt{\det(v+t)}}\mu(dv)=\frac{1}{2^{1+\d n}}\frac{t^{-1}}{\sqrt{\det t}}.\end{equation}

\end{enumerate}

\vspace{4mm}\noindent\textbf{Proof.} We have 
 \begin{equation}
 \label{fhat}
 \hat{f}(s)=\int_{\R^n}e^{i\<s,x\>}f(x)dx=\E(e^{i\<V^{1/2}Z,s\>})=\E(e^{-\d s^*Vs}).
 \end{equation}
 Now using Plancherel Theorem and Lemma 5.1, we prove part 1. as follows:
$$\int_{\R^n}f^2(x)dx=\frac{1}{(2\pi)^n}\int_{\R^n}\hat{f}(s)^2ds
=\frac{1}{(2\pi)^n}\int_{\R^n}\E(e^{-\d s^*(V+V_1)s})ds=\frac{1}{(2\pi)^{n/2}}\E\left(\frac {1}{\det\sqrt{V+V_1}}\right).$$

To  prove part 2, we use Plancherel theorem again for the function 
$$g(x)=f(x)-\frac{e^{-\frac{x^*t^{-1}x}{2}}}{(2\pi)^{n/2}(\det t)^{1/2}},$$
and obtain
$$I(t)=\frac{1}{(2\pi)^{n}}\int_{\R^n}\left[\hat{f}(s)-h_t(s)\right]^2ds$$where $h_t=e^{-\d s^*ts}.$ From Lemma 4.1 applied to $A=2t$ we have $\|h_t\|^2=\pi^{n/2}/\sqrt{\det t}.$
Expanding the square in $I(t)$  we obtain
$$(2\pi)^nI(t)-||\hat{f}||^2=\frac{(\pi)^{n/2}}{\sqrt{\det t}}-2\langle \hat{f}, h_t\rangle:=I_1(t),$$
where $h_t=e^{-\d s^*ts}.$
We now want to show that the minimum of $I_1(t)$ is reached at some $t_0\in \mathcal{P}.$ 

We show that 
$$K_1=\{y\in \mathcal{P}; I_1(y^{-1})\leq 0\}$$ is  non empty and compact. Writing 
$$I_2(y)=\<\hat{f},h_{y^{-1}}\>\frac{1}{(2\pi)^{n/2}\sqrt{\det y}},$$ 
we see that $y\in K_1$, i.e.  $I_1(y^{-1})\leq 0$ if and only if $\frac{1}{2^{1+\d n}}\leq I_2(y).$
From \eqref{fhat}, the definition of $h_t(s)$ and Lemma 4.1, we have that
\begin{eqnarray*}
I_2(y)&=&\frac{1}{(2\pi)^{n/2}\sqrt{\det y}}\int_{\R^n}\E(e^{-\frac{s^*Vs}{2}})e^{-\frac{s^*y^{-1}s}{2}}ds=\frac{1}{(2\pi)^{n/2}\sqrt{\det y}}\E\Big(\int_{\R^n}
e^{\frac{s^*(V+y^{-1})s}{2}}ds\Big)\\
&=&\frac{1}{\sqrt{\det y}}\E\Big(\frac{1}{\sqrt{\det (V+y^{-1})}}\Big)=\int_{\mathcal{P}}\frac{\mu(dv)}{\sqrt{\det(I_n+vy)}}.
\end{eqnarray*}
For $0<C\leq 1$
let us show that 
$$K_2=\{y\in \mathcal{P}; I_2(y)\geq C\}$$  is compact. Note that $K_1=K_2$ for  $C=1/2^{1+\d n}.$ Since $I_2$ is continuous, $K_2$ is closed. The set $K_2$ is  not empty since 
$I_2(y)\geq 1.$ Let us prove that $K_2$ is bounded. Recall $\|y\|=(\tr y^2)^{1/2}.$ Suppose that $y^{(k)}\in K_2$ is such that $\|y^{(k)}\|\to_{k\to \infty} \infty$ and let us show that for such a $y^{(k)},\;I_2(y^{(k)})\to 0,$ which is a contradiction.

 Indeed, $\tr (vy^{(k)})\to_{k\to \infty} \infty$ if $v\in \mathcal{P}$ . To see this, assume that $v=\mathrm{diag}(v_1,\ldots,v_n)$. Then
 \begin{eqnarray*}\tr (vy^{(k)})&=&v_1y_{11}^{(k)}+\cdots+v_ny_{nn}^{(k)}\\&\geq& \tr (y^{(k)})\times\min_i v_i\geq \|y^{(k)}\|\times\min_i v_i\to_{k\to \infty} \infty,\end{eqnarray*}
 where the last inequality is due to the fact that if  $\lambda_1,\ldots, \lambda_n$ are positive, then $\sqrt{\lambda_1^2+\ldots,+\lambda_n^2}\leq \lambda_1+\ldots+\lambda_n$. 
 Moreover, if $(\lambda_1.\ldots,\lambda_n)$ are the eigenvalues of $vy^{(k)}$,
$$\det (I_n+vy^{(k)})=(1+\lambda_1)\ldots(1+\lambda_n)\geq 1+\lambda_1+\cdots+\lambda_n=1+\tr (vy^{(k)})\to_{k\to \infty} \infty$$ By dominated convergence, it follows that $I_2(y^{(k)})\to_{k\to \infty}0$ and this proves that $K_2$ is bounded. We have therefore shown that $K_1$ is  compact. This proves that the minimum of $I_1(t)$ and thus of $I(t)$ is reached at some point $t_0$ of $\mathcal{P}.$

The last task is to show that $t_0$ is a solution of  equation $\eqref{DIFFn}.$ Since $I(t)$ is differentiable and reaches its minimum on the open set $\mathcal{P}$, the differential of $I(t)$ must cancel at $t_0$. The differential of $I$ is the following linear form on $\mathcal{S}$
$$h\in \mathcal{S} \mapsto I'(t)(h)=\frac{1}{(2\pi)^{n/2}}\int_{\R^n}\left[\hat{f}(s)-e^{-\d s^*ts}\right]e^{-\d s^*ts}s^*hsds.$$  The equality $I'(t)=0$ is equivalent to
$$\int_{\R^n}\hat{f}(s)e^{-\d s^*ts}ss^*ds=\int_{\R^n}e^{- s^*ts}ss^*ds.$$ Using the second formula in Lemma 4.1 and the fact that $\hat{f}(s)=\E(e^{-\d s^*Vs})$, we obtain

$$\int_{\mathcal{P}}\frac{(v+t)^{-1}}{\sqrt{\det(v+t)}}\mu(dv)=\frac{(2t)^{-1}}{\sqrt{\det (2t)}}=\frac{1}{2^{1+\d n }}\frac{t^{-1}}{\sqrt{\det t}},$$
which proves \eqref{DIFFn}. 

\vspace{4mm}\noindent \textbf{Remarks.} \begin{enumerate}\item We note that \eqref{DIFFn} can also be written in terms of $y=t^{-1}$ as
$$\int_{\mathcal P} \frac{(1+vy)^{-1}}{\sqrt{\det(1+vy)}}\mu(dv)=\frac{1}{2^{1+\frac{n}{2}}}I_n.$$

\item While it is highly probable that the value $t_0$ at which $I(t)$ reaches its minimum is unique, it is difficult to show for $ n\geq 2$ that   equation \eqref{DIFFn} has a unique solution: there is no reason to think that the function $ t\mapsto I(t) $ is convex. However a case of uniqueness is proved in Proposition 5.3 below. \end{enumerate}

\subsection{Best approximation for a scalar mixture.}

\vspace{4mm}\noindent\textbf{Proposition 5.3.} 
Let $\nu(d\lambda)$ be a probability on $(0,\infty)$ such that $$\E((\Lambda+\Lambda_1)^{-n/2})<\infty$$ where $\Lambda$ and $\Lambda_1$ are independent with distribution $\nu$, and let $\mu$ be the distribution of $V=\Lambda I_n.$ Then $t\mapsto I(t)$ defined in \eqref{ioft} reaches its minimum at a unique point $t_0$. Furthermore $t_0$ is a multiple of $I_n.$

\vspace{4mm}\noindent\textbf{Proof.} From Theorem 4.2, $I$ reaches its minimum at least at one point $t_0\in \mathcal {P}.$ Without loss of generality by choosing a suitable orthonormal basis of $\R^n$, we can assume that $t_0=\mathrm{diag}(\lambda^0_1,\ldots,\lambda^0_n).$ We are going to show that  $\lambda^0_1=\ldots=\lambda^0_n.$ Consider the restriction $I^*$ of $I$ to the set of  diagonal matrices with positive entries, namely
$$I^*(t_1,\ldots,t_n)=I^*(\mathrm{diag}(t_1,\ldots,t_n)).$$ Of course $(t_1,\ldots,t_n)\mapsto I^*(t_1,\ldots,t_n)$ reaches its minimum on $(\lambda^0_1,\ldots,\lambda^0_n).$ 
By a computation which imitates the proof of Theorem 4.2 we consider
\begin{eqnarray*}
I_1^*(t_1,\ldots,t_n)&=&(2\pi)^nI^*(t_1,\ldots,t_n)-\|\hat{f}\|^2\\
&=&\frac{\pi^{n/2}}{\sqrt{t_1\ldots t_n}}-2\int_0^{\infty}\frac{\nu(d\lambda)}{\prod_{i=1}^n(t_i+\lambda)^{1/2}}
\end{eqnarray*}
Since  $I_1^*(t_1,\ldots,t_n)$ reaches its minimum at $t_0$, its gradient is zero at $(\lambda^0_1,\ldots,\lambda^0_n)$. We have
$$\frac{\partial}{\partial t_j}I_1^*(t_1,\ldots,t_n)=-\frac{\pi^{n/2}}{2t_j\sqrt{t_1\ldots t_n}}+\int_0^{\infty}\frac{\nu(d\lambda)}{(t_j+\lambda)\prod_{i=1}^n(t_i+\lambda)^{1/2}}$$
and  as a consequence,  for all $j=1,\ldots,n$
\begin{equation}\label{ASTUCE}
\int_0^{\infty}\frac{\lambda^0_j}{\lambda^0_j+\lambda}\times \frac{\nu(d\lambda)}{\prod_{i=1}^n(\lambda^0_i+\lambda)^{1/2}}=\frac{\pi^{n/2}}{2\sqrt{\lambda^0_1\ldots \lambda^0_n}}.
\end{equation} 
The important point of \eqref{ASTUCE} is the fact that the right hand side does not depend on $j.$ Suppose now that there exists $j_1$ and $j_2$ such that $\lambda^0_{j_1}<\lambda^0_{j_2}.$ This implies that for all $\lambda>0$ we have 
$$\frac{\lambda^0_{j_1}}{\lambda^0_{j_1}+\lambda}<\frac{\lambda^0_{j_2}}{\lambda^0_{j_2}+\lambda}$$ and the left hand sides of  \eqref{ASTUCE} cannot be equal for $j=j_1$ and $j=j_2.$ As a consequence $t_0=\lambda^0I_n$ for some $\lambda^0>0.$

To see that $\lambda^0$ is unique, we imitate the proof of Theorem 3.1. We omit the details here.    \hfill $\square$

We will finish by giving an example of a scalar Gaussian mixture, actually built on the univariate Kolmogorov-Smirnov measure \eqref{KS} with density
$$k_1(\lambda)=\sum_{n=1}^{+\infty} (-1)^{n+1}n^2e^{-\frac{n^2\lambda}{2}}{\bf 1}_{(0,+\infty)}(\lambda).$$

\vspace{4mm}\noindent \textbf{Example 5.} Let us verify first that 
$$g_n(x)=C_n\frac{e^{||x||}}{(1+e^{||x||})^2},$$
where $C_n$ is the normalizing constant, is a density in $\R^n$. Indeed, using polar coordinates in $\R^n$ with $r=||x||$, we have $\frac{1}{C_n}=S_{n-1}J(n-1)$ where $S_{n-1}=n\pi^{n/2}/\Gamma(1+\frac{n}{2})$ is the area of the unit sphere in $\R^n$ and where $$J(t)=\int_0^{+\infty}\frac{e^{-r} r^{t}}{(1+e^{-r})^2}dr.$$ Of  course $J(0)=1/2$ and by integration by part $J(1)=\log 2.$ For $t>1$ we have 
\begin{eqnarray*}
J(t)&=&\sum_{k=1}^{\infty}(-1)^{k-1}k\int_0^{\infty}e^{-kr}r^tdr=\Gamma(t+1)\sum_{k=1}^{\infty}\frac{(-1)^{k-1}}{k^t}\\
&=&\Gamma(t+1)
(1-2^{{1-t}})\zeta(t).
\end{eqnarray*}
where $\zeta(t)=\sum_{k=1}^{\infty}\frac{1}{k^t}$ is the Riemann function and the last equality is a well-known formula. Thus for instance $$C_1=1, C_2=1/(2\pi \log 2),\ C_3=3/(2\pi^3).$$ Next, writing
$$k_n(\lambda)=C_n(2\pi\lambda)^{\frac{n-1}{2}}k_1(\lambda){\bf 1}_{(0, +\infty)}(\lambda)$$
let us show that $k_n$ is a density such that
\begin{equation}
\label{kn}
\int_0^{+\infty}\frac{e^{-\frac{||x||^2}{2\lambda}}}{(2\pi\lambda)^{n/2}}k_n(\lambda)d\lambda=g_n(x).
\end{equation}
This means, of course, that $g_n$ is a scale mixture of multivariate normal $N(0, \lambda I_n)$ distributions. We have 
\begin{eqnarray*}
1=\int_{\R^n}g_n(x)dx&=&C_n\int_{\R^n} \int_0^{+\infty}\frac{e^{-\frac{||x||^2}{2\lambda}}}{(2\pi\lambda)^{n/2}}(2\pi\lambda)^{(n-1)/2}k_1(\lambda) d \lambda dx\\
&=&C_n\int_0^{+\infty}(2\pi\lambda)^{(n-1)/2}k_1(\lambda) \left(\int_{\R^n}\frac{e^{-\frac{||x||^2}{2\lambda}}}{(2\pi\lambda)^{n/2}}dx\right) d\lambda\\
&=&C_n\int_0^{+\infty}(2\pi\lambda)^{(n-1)/2} k_1(\lambda) d \lambda.
\end{eqnarray*}

\section{References}
\vspace{4mm}\noindent\textsc{Feller, W.} (1966) \textit{An Introduction to Probability Theory and its Applications}, Vol. 2, Wiley, New York.

\vspace{4mm}\noindent\textsc{Gneiting, T.} (1997) Normal scale mixture and probability densities. \textit{J. Stat. Comput Simul.} \textbf{59}, 375-384.

\vspace{4mm}\noindent\textsc{Hardy, G. H.,  Wright, E. M.} (1938) \textit{An Introduction to the Theory of Numbers}, Heath-Brown, D. R., Silverman, J. H., eds.,  (Sixth ed. 2008), Oxford University Press.

\vspace{4mm}\noindent\textsc{Kolmogorov, A. N.} (1933) Sulla determinazi\'one emp\`irica di una l\'egge di distribuzi\'one. \textit{G. Inst. Ital. Attuari} \textbf{4}, 83-91.

\vspace{4mm}\noindent\textsc{Letac, G., Massam, H. and Mohammadi, R.} (2017) The ratio of normalizing constants for Bayesian Gaussian model selection,  revision oct. 12th  2018, arXiv 1706. 04416.

\vspace{4mm}\noindent\textsc{Monahan, J. F. and Stefanski, L. A.} (1992)  Normal Scale
Mixture Approximations to $F^*(z)$ and Computation of the
Logistic-Normal Integral, in {\it Handbook of the Logistic
Distribution}, N. Balakrishnan, Ed., Marcel Dekker, New York.

\vspace{4mm}\noindent\textsc{Palmer, J. A., Kreutz-Delgado, K. and Makeig, S.} (2011) Dependency models based on generalized Gaussian scale mixtures. \textit{DRAFT UCSD-SCCN v1.0}, Sept 7.

\vspace{4mm}\noindent\textsc{Stefanski, L. A.} (1991) A Normal Scale Mixture
Representation of the Logistic Distribution, {\it Statistics \&
Probability Letters} {\bf 11}, 69--70. 

\vspace{4mm}\noindent\textsc{West, M.} (1987) On scale mixture of normal distributions, \textit{Biometrika} \textbf{74}, 3, 646-648.

\end{document}